\documentclass[11pt]{article}
\usepackage{amsmath}
\usepackage{amsfonts}
\usepackage{amssymb}
\usepackage{epsfig}
\usepackage{latexsym}

\newtheorem{theorem}{\bf Theorem}[section]
\newtheorem{lemma}[theorem]{\bf Lemma}

\def\sgn{\mbox{sgn}}

\begin{document}

\title{Two-Dimensional Adaptive Fourier Decomposition}
\author{Tao~Qian \thanks{Department of Mathematics, University of Macau, Macao (Via Hong
Kong). Email: fsttq@umac.mo. Telephone: +853 83978547. Fax: +853 28838314. The work was supported by
Multi-Year Research Grant (MYRG) MYRG116(Y1-L3)-FST13-QT, Macao Government FDCT/056/2010/A3 and FDCT 098/2012/A3}}
\maketitle

\begin{abstract} One-dimensional adaptive Fourier decomposition, abbreviated as 1-D AFD, or AFD, is an adaptive
 representation of a physically realizable signal into a linear combination of parameterized Szeg\"o and higher order Szeg\"o kernels of the context. In the present paper we study multi-dimensional AFDs based on multivariate complex Hardy spaces theory.  We proceed with two approaches of which one uses Product-TM Systems; and the other uses Product-Szeg\"o Dictionaries. With the Product-TM Systems approach we prove that at each selection of a pair of parameters the maximal energy may be attained, and, accordingly, we prove the convergence. With the Product-Szeg\"o dictionary approach we show that Pure Greedy Algorithm is applicable. We next introduce a new type of greedy algorithm, called Pre-Orthogonal Greedy Algorithm (P-OGA). We prove its convergence and convergence rate estimation, allowing a weak type version of P-OGA as well. The convergence rate estimation of the proposed P-OGA evidences its advantage over Orthogonal Greedy Algorithm (OGA). In the last part we analyze P-OGA in depth and introduce the concept P-OGA-Induced Complete Dictionary, abbreviated as Complete Dictionary . We show that with the Complete Dictionary P-OGA is applicable to the Hardy $H^2$ space on $2$-torus. \end{abstract}


\vspace *{0.5mm}

\section {Preparation}

We will give a brief introduction to the related background knowledge as well as the existing 1-D AFD theory. Denote by ${\bf D}$ the unit disc and ${\bf C}$ the complex plane. The present paper concentrates in the unit disc context. There is a parallel theory for the upper-half plane context (\cite{QWa}, \cite{QWang}). The space $L^2({\bf \partial D})$ can be expressed as the direct sum of the two relevant boundary Hardy spaces, viz.,
\[ L^2({\bf \partial D})=H^2_+({\bf \partial D})\bigoplus H^2_-({\bf \partial D}),\]
where $H^2_+ ({\bf \partial D})$ and $H^2_- ({\bf \partial D})$ consist of, respectively, the non-tangential boundary limits of the complex Hardy $H^2$-functions inside and outside the unit disc. The mentioned complex Hardy spaces of holomorphic functions inside and outside the unit disc are, respectively, denoted by $H^2_\pm ({\bf D}).$ The non-tangential boundary limit mappings from $H^2_\pm ({\bf D})$ to their boundary limit spaces $H^2_\pm ({\bf \partial D})$ are isometric isomorphisms. The closed subspaces $H^2_\pm ({\bf \partial D})$ of $L^2({\bf \partial D}))$ are, in fact, the collections of the functions of the forms $(f\pm i Hf)/2\pm c_0/2, c_0\in {\bf C},$ respectively, where $H$ is the Hilbert transformation of the context (see below). Signals of the above forms are called analytic signals (\cite{Gab}).

The Hilbert transformation operator for a simply-connected region $\Omega$ is related to the Plemelj formula of the context  (\cite{Be}). Let $f$ be a $\lq\lq$good" holomorphic function in $\Omega$ such that it has non-tangential boundary limits almost everywhere on the boundary $\partial \Omega.$ By denoting the boundary limit function as $\tilde{f}=u+iv,$ where $u$ and $v$ are real-valued, we call $v$ \emph{the Hilbert transform of} $u.$ Note that both $u$ and $v$ are defined on $\partial \Omega.$  It is a fundamental result that functions in the complex Hardy spaces have non-tangential boundary limits almost everywhere on the boundary. The Hilbert transforms of the square integrable functions on the unit disc are given by
\[ Hu(t)= \frac{1}{2\pi}\lim_{\epsilon\to 0}\int_{|t-s|>\epsilon} \cot \left(\frac{t-s}{2}\right)u(e^{is})ds\]
  (\cite{Ga}). The Hilbert transformation has a Fourier multiplier representation:
 \[ Hu (t)=\sum_{n=-\infty}^\infty (-i \sgn (n)) c_n e^{int}, \quad u(t)=\sum_{n=-\infty}^\infty c_n e^{int}, \quad \sum_{n=-\infty}^\infty |c_n|^2<\infty,\]
where $\sgn (\xi)$ is the signum function, being of the value $+1$ or $-1,$ respectively, for $\xi>0$ or $\xi<0;$  and $\sgn (0)=0.$

 From now on we assume that functions to be studied in $L^2({\bf \partial D})$ are real-valued.  Under such assumption while both $f$ and $Hf$ are real-valued. Denote $f^\pm=\frac{1}{2}\left(f\pm iHf\right)\pm c_0/2.$ Due to the relation $c_{-n}=\overline{c_n}$ we have
 \begin{eqnarray}\label{real valued 1} f=f^+ + f^-, \qquad f=2{\rm Re}f^+-c_0.\end{eqnarray}
 We note that for any $f\in L^2({\bf \partial D}),$ not necessarily real-valued, the Plemelj formula inside the disc is
 \[ \lim_{r\to 1-}\frac{1}{2\pi i}\int_{\bf \partial D}\frac{f(\zeta)}{\zeta -re^{it}}d\zeta = \frac{1}{2}(f(e^{it})
 +i Hf (e^{it}))+\frac{c_0}{2}, \qquad {\rm a.e.},\]
 and, in the $L^2$-convergence sense,
 \[ f^+(e^{it})=\sum_{k=0}^\infty c_ke^{ikt}, \qquad f^-(e^{it})=\sum_{k=-\infty}^{-1} c_ke^{ikt}.\]

 For a real-valued function $f$ the second relation in (\ref{real valued 1}) shows that a series expansion of $f^+$ gives rise to a series expansion of $f.$
 AFD is an adaptive decomposition (expansion) for $f^+$ into a linear combination of parameterized Szeg\"o and higher order Szeg\"o kernels of the context. On the disc the $L^2$-normalized Szeg\"o kernels alone constitute a dictionary, called \emph{Szeg\"o Dictionary},
 \[ {\cal D}=\{e_a\}_{a\in {\bf D}}, \quad e_a(z)=\frac{\sqrt{1-|a|^2}}{1-\overline{a}z}.\]
 A Szeg\"o kernel has the remarkable property that it is the reproducing kernel of the underlying Hardy space: Under the inner product
\[ \langle f, g \rangle = \frac{1}{2\pi}\int_0^{2\pi} f(e^{it})\overline{g}(e^{it})dt,\] there holds
\begin{eqnarray}\label{reproducing}\langle f, e_a \rangle =\sqrt{1-|a|^2}f(a).\end{eqnarray}
 We note that AFD not only gives a series in Szeg\"o and higher order Szeg\"o kernels that fast converges to the targeted function (\cite{QWang}), but also gives rise to a positive frequency representation of the function with considerable stability (\cite{Pi}, \cite{QWa}). The recent studies of positive instantaneous frequency and related signal decomposition were motivated by related work in signal analysis area, including those of Picinbono, Cohen and Huang (\cite{Pi}, \cite{Co}, \cite{Huang}). We note that in spite of the desires and efforts, EMD of Huang et al does not generate positive analytical instantaneous frequency, nor what they called HHT (\cite{SV}). Studies show that theoretical development of this theme has to involve complex Hardy spaces, and, especially, the inner and outer functions theory of Nevanlinna in complex analysis.

In rational approximation of one complex variable one cannot avoid the so called
 rational orthogonal systems, or alternatively, Takenaka-Malmquist (TM) systems
(\cite{Wa}, \cite{Bu}). The functions in a TM system are obtained, in fact, from the Gram-Schmidt (G-S) orthogonalization process applied to the \emph{partial fractions} $E_k,$ essentially the parameterized Szeg\"o and higher order Szeg\"o kernels, defined as follows. Let $\{a_1,...,a_n\}$ be an $n$-tuple of complex numbers in ${\bf D}.$ We say that an entry $a_k$ of the $n$-tuple $\{a_1,\ldots,a_n\}$
\emph{has the multiplicity} $m_k$, if there exist exactly $m_k$ entries
$a_{n_1}, \ldots, a_{n_{m_k}}$ with $1\leq n_1<\cdots<n_{m_k}=k$ such that
$a_{n_1}=\cdots =a_{n_{m_k}}=a_k$.
Given an $n$-tuple $\{a_1,\ldots,a_n\},$ we define
\begin{equation} \label{def.ps}
E_k(z) = E_{\{a_1,\ldots,a_k\}}(z)
:= \left\{\begin{array}{lll} \displaystyle \frac{1}{(1-\overline{a}_kz)^{m_k}}
&\text{if} & a_k\not=0 \\
z^{m_k-1} & \text{if} & a_k=0 \\
\end{array}\right.
\end{equation}
where $m_k$ is the multiplicity of $a_k, k=1,...,n.$ We say that $E_k$ is a
\emph{higher order Szeg\"o kernel} if $m_k>1;$ and otherwise a \emph{Szeg\"o kernel}.
The system $\{E_1,\ E_2,\ \ldots,\ E_n\}$ is called the \emph{partial fraction system generated by} $\{a_1,\ldots,a_n\}.$  It is usually not orthogonal.
 A TM system, as denoted by $\{B_k\}_{k=1}^n$ in the sequel, is the result of the G-S orthogonalization process applied to a partial fraction system $\{E_k\}_{k=1}^n$ (\cite{Bu}, \cite{QWe}).  In the sequel we will call the dictionary consisting of all the normalized partial fractions, viz., the normalized Szeg\"o and higher order Szeg\"o kernels given in (\ref{def.ps}), the \emph{Complete Szeg\"o Dictionary}, denoted by $\tilde{\cal D}$ (see \S 4).

 In this article we generalize the theory in the unit disc to poly-discs. The $2$-disc theory, in particular, like multiple Fourier series, has direct applications to image analysis and image processing.

    In the unit disc context, a TM system is an infinite sequence of parameterized rational functions
\[ B_k(z)=\frac{\sqrt{1-|a_k|^2}}{1-\overline{a}_kz}\prod_{l=1}^{k-1}\frac{z-a_l}{1-\overline{a}_lz}, \ k=1,2,...,\]
where $a_1,...,a_k,...$ are any complex numbers in the open unit disc. The multiple product part of each $B_k$ is a Blaschke product with $k-1$ zeros, being the product of the explicitly given $(k-1)$ M\"obius transforms, while the rest part is a normalized Szeg\"o kernel, being an element of the dictionary ${\cal D},$ where we treat the complex number $a$ as a parameter. Such systems have been well studied. In particular, when all the $a_k$'s are identical with zero,  $\{B_k\}$ reduces to a half of the Fourier system, viz., $\{z^{k-1}\}_{k=1}^\infty,$ a basis of the Hardy spaces. It is known that for general parameters $a_k$'s a system $\{B_k\}$ is a basis of the Hardy space $H^p$ on the disc, $1\leq p\leq \infty,$ if and only if the hyperbolic non-separability condition is satisfied, viz.,
\begin{eqnarray}\label{non-separate condition} \sum_{k=1}^\infty (1-|a_k|)=\infty\end{eqnarray}
(\cite{Bu}). For $p=2$ the TM system is an orthonormal basis of $H^2({\bf D})$ if and only if the condition (\ref{non-separate condition}) is met.

It is noted that writing in the form $B_k(e^{it})=\rho_k(t)e^{i(\psi_k(t)+\phi_{k-1}(t))},$ where $e_{a_k}(e^{it})=\rho_k(t)e^{i\psi_k(t)}, \rho_k(t)\ge 0, e^{i\phi_{k-1}(t)}=\prod_{l=1}^{k-1}\frac{e^{it}-a_l}{1-\overline{a}_le^{it}},$ we have $\phi_{k-1}'(t)\ge 0$ and $1+\psi_k'(t)>0$ for all $t\in [0,2\pi ], k=1,2,..., \phi'_0=0.$  In particular, if $a_1=0,$ then $\psi'_k(t)+\phi'_{k-1}(t)\ge 0$ for all $t$ and $k=1,2,...$ (\cite{QWa}).  A representation into a linear combination of such $B_k$'s is a positive instantaneous frequency decomposition, or alternatively, a mono-component (Hardy space functions of positive analytic phase derivative almost everywhere) decomposition (\cite{Pi}, \cite{Co}, \cite{QWa}, \cite{QWYZ}, \cite{Tan-Shen-Yang}, \cite{YuZ}). In fact, it was the seeking for such decompositions that motivated 1-D AFD.

To make easy understanding to 2-D AFD delivered in \S 2 and also for the self-containing purpose, we now give an exposition for the existing 1-D AFD
(\cite{QWa}, \cite{QWang}). Let $f$ belong to the Hardy space $H^2({\bf D}).$ Set $f_1=f.$
For any $a_1\in {\bf D},$  we have the identity
 \begin{eqnarray}\label{maximum sifting}  f(z)=\langle f_1, e_{a_1}\rangle e_{a_1}(z)+f_2(z)\frac{z-a_1}{1-\overline{a}_1z},\end{eqnarray}
with
\[f_2(z)=\frac{f_1(z)-\langle f_1, e_{a_1}\rangle e_{a_1}(z)}{\frac{z-a_1}{1-\overline{a}_1z}}.\]
Due to the reproducing kernel property of $e_{a_1}$ in $H^2({\bf D})$ we have
\[ \langle f, e_{a_1} \rangle =\sqrt{1-|a_1|^2}f(a_1), \quad {\rm and\ hence,} \quad f_1(a_1)-\langle f_1, e_{a_1}\rangle e_{a_1}(a_1)=0.\]
The last assertion implies that $f_2\in H^2({\bf D})$ that enables the recursive process in the sequel.
We call the transformation from $f_1$ to $f_2$ the \emph{generalized backward shift via }$a_1;$ and $f_2,$ the \emph{reduced remainder}, being \emph{the generalized backward shift transform of} $f_1$ \emph{via} $a_1$. Such terminology was motivated by the classical backward shift operator
\[ S(f)(z)=\sum_{k=0}^\infty c_{k+1}z^k  = \frac{f(z)-f(0)}{z},\]
where we assume $f(z)=\sum_{k=0}^\infty c_kz^k.$
Noticing that $f(0)=\langle f, e_0 \rangle e_0(z),$  $S$ is identical with the just defined
 generalized backward shift operator via $0.$

Due to the orthogonality between the two terms on the right hand side of (\ref{maximum sifting})
 and the unimodular property of M\"obius transforms on the boundary, we have
\[ \|f\|^2=\|\langle f_1, e_{a_1} \rangle e_{a_1}\|^2 + \|f_2\frac{(\cdot ) -a_1}{1-\overline{a}_1(\cdot )}\|^2=|\langle f_1, e_{a_1} \rangle |^2+\| f_2\|^2.\]
  We are to select $a_1\in {\bf D}$ that gives the term $\langle f_1, e_{a_1} \rangle e_{a_1}(z)$ the maximal energy
   out of  $\| f\|^2.$ Due to the reproducing kernel property of $e_a$ we have
\begin{eqnarray}\label{reason1} | \langle f_1, e_{a_1}\rangle |^2=(1-|a_1|^2)|f_1(a_1)|^2,\end{eqnarray}
and, hence, we are reduced to find
$a_1\in {\bf D}$ such that
\[ a_1=\arg \max \{ (1-|a|^2)|f_1(a)|^2\ : \ a\in {\bf D} \}.\]
The existence of such maximal selection is evident (\cite{QWa}, it can also be referred to the proofs of Theorem \ref{maximal} and Theorem \ref{Maximal Szego} below), and called \emph{Maximal Selection Principle}.  Having selected such $a_1$ we repeat the same procedure to $f_2,$ and so on.  After consecutive
$n$ steps, we obtain
\[ f(z)=\sum_{k=1}^n \langle f_k, e_{a_k}\rangle B_k(z) + f_{n+1}\prod_{k=1}^n\frac{z-a_k}{1-\overline{a}_kz},\]
where for $k=1,...,n,$
\begin{eqnarray}\label{reason2} a_k=\arg \max \{ (1-|a|^2)|f_k(a)|^2\ : \ a\in {\bf D}\},\end{eqnarray}
and, for $k=2,...,n+1,$
\[f_k(z)=\frac{f_{k-1}(z)-\langle f_{k-1}, e_{a_{k-1}}\rangle e_{a_{k-1}}(z)}{\frac{z-a_{k-1}}{1-\overline{a}_{k-1}z}}.\]
Due to the orthogonality, we have
\[ \|f-\sum_{k=1}^n \langle f_k, e_{a_k}\rangle B_k(z) \|^2=\| f\|^2-\sum_{k=1}^n |\langle f_k, e_{a_k}\rangle|^2=\| f_{k+1}\|^2.\]
It can be shown that
\[ \lim_{n\to \infty} \| f_{k+1} \|=0\] (\cite{QWa}), and thus
\begin{eqnarray}\label{AFD}  f(z)=\sum_{k=1}^\infty \langle f_k, e_{a_k}\rangle B_k(z).\end{eqnarray}
The decomposition of $f$ given by (\ref{AFD}) is called Adaptive Fourier Decomposition (AFD) of $f.$

The following relations are noted:
\begin{eqnarray}\label{reduce} \langle f_k, e_{a_k}\rangle =\langle g_k,B_k\rangle =\langle f,B_k\rangle ,\end{eqnarray}
where $g_k$ is the \emph{orthogonal standard remainder} defined through
\begin{eqnarray}\label{standard} f=\sum_{i=1}^{k-1} \langle f, B_i\rangle B_i(z)+g_k(z).\end{eqnarray}
We also cite the useful relations
\begin{eqnarray}\label{relations}     g_k(z)=f_k(z)\prod_{l=1}^{k-1}\frac{z-a_l}{1-\overline{a}_lz},\quad
{\rm where} \quad f_k=S_{a_{k-1}}
\cdots S_{a_1}f(z).\end{eqnarray}

To deal with the convergence rate issue we define a particular subclass
of functions (\cite{QWang}):
\begin{equation}
H^{2}(\mathcal {D},M):=\{f\in H^2(D):
f=\sum_{k=1}^{\infty}c_{k}e_k, e_{k}\in\mathcal {D}, \
\sum_{k=1}^{\infty}|c_{k}|\leq M\}, \quad 0<M<\infty.
\end{equation}
We have (\cite{QWang})

\begin{theorem}
Let $\mathcal {D}$ be the dictionary of the normalized Szeg\"o kernels
of $H^2(D)$. Then for each $f\in H^{2}(\mathcal {D},M)$,
decomposed by Adaptive Fourier Decomposition, we have
\begin{eqnarray}
\parallel g_k \parallel \leq \frac{M}{\sqrt{k}}. \nonumber
\end{eqnarray}\end{theorem}

We note that functions in the class $H^{2}(\mathcal {D},M)$ may not be smooth on the unit circle. The above estimation in the energy sense reflects the tolerance of AFD with non-smoothness.\\

\noindent {\bf Remark 1} There are variations of 1-D AFD (also called Core AFD in later references as it is the construction block of the other AFDs based on maximal selections) of which we mention Unwinding AFD (\cite{Q}), Cyclic AFD (\cite{Qcyclic}) and  Higher-Order-Szeg\"o-Kernel AFD (\cite{QWj}). Each of the mentioned AFD variations has its own merits. In particular, Unwinding and Higher-Order-Szeg\"o-Kernel AFDs are designed for decompose signals of high frequencies. Cyclic AFD offers a conditional solution for the open problem of finding a rational Hardy space function whose degree does not exceed a pre-described integer $n$ that best approximates a given Hardy space function (\cite{Qcyclic}).\\

\noindent {\bf Remark 2}
1-D AFD and its variations have found significant applications to system identification and signal analysis (\cite{MQ}, \cite{MQW}, \cite{MQ2}, \cite{QChenTan}).\\

\noindent {\bf Remark 3} Any complete TM system with $a_1=0$ gives rise to signal decompositions of positive frequencies.  To sufficiently characterize a signal it is desirable to find the most suitable TM systems for the given signal. The suitability may be measured by the corresponding convergence speed. It was the effort of gaining fast convergence that made AFD to share the same idea as greedy algorithm.\\

\noindent{\bf Remark 4} 1-D AFD, as a matter of fact, is not the same as any existing greedy algorithm. In the last section of this article we will introduce a new greedy algorithm called \emph{Pre-Orthogonal Greedy Algorithm} (\emph{P-OGA}) and will show that P-OGA is identical with AFD in the unit disc context.  Among the existing greedy algorithms, including the weak type ones, the one that is most close to P-OGA is Orthogonal Greedy Algorithm (OGA). In AFD, at the step $n,$ a parameter value $a_n$ can be repeatedly selected with respect to the reduced remainder $f_n$ in order to have the maximal energy gain.  On the other hand,
OGA does not allow repeated selections of a parameter with respect to the orthogonal standard remainder $g_n.$ In the unit disc context OGA generates orthogonal projections of the orthogonal standard remainder $g_n$ into linear span of
\[ 1, \ \frac{1}{1-\overline{a}_1z},\ ..., \ \frac{1}{1-\overline{a}_nz}, \quad n=1,2,...,\]
where all the $a_k$'s are distinct.  In contrast, 1-D AFD, in accordance with (\ref{def.ps}),  gives rise to orthogonal projections into linear spans of the partial fraction systems
\begin{eqnarray}\label{higher} 1, \ ...\ , z^{m_0-1}, \frac{1}{1-\overline{a}_1z}, \ ...,\ \frac{1}{(1-\overline{a}_1z)^{m_1}}, \ ..., \ \frac{1}{1-\overline{a}_nz}, \ ...,\ \frac{1}{(1-\overline{a}_nz)^{m_n}}, \quad n=1,2,...,\end{eqnarray}
where all the $a_n$ are distinct and $m_n$ are the respective multiples.
The latter is with the full strength of the related partial fractions and thus the decomposition converges faster (also see \S 2 and \S 3). Not only having a fast converging positive frequency representation, AFD through its backward shift process                               automatically generates an orthogonal expansion without involving the G-S orthogonalization process.\\

\noindent {\bf Remark 5} The Maximal Selection Principle may produce a sequence $a_1,...,a_n,...$ that does not satisfies (\ref{non-separate condition}).  Such case corresponds to a remarkable decomposition of the Hardy space. It is noted that the case is exactly when a Blaschke product $\phi (z)$ is definable with $a_1,...,a_k,...$ being all its zeros including the multiples (\cite{Ga}).  In such case we have the space  decomposition
\begin{eqnarray}\label{above relation} H^2=\overline{{\rm span}\{B_k\}}\oplus \phi H^2,\end{eqnarray}
where the closed set $\overline{{\rm span}\{B_k\}}$ is a \emph{backward shift invariant subspace} and $\phi H^2$ a \emph{shift invariant subspace} of the $H^2$ space. In such case $f$ belongs to the backward shift invariant space. Backward shift invariant subspaces have significant applications to phase and amplitude retrieval problems and solutions of the Bedrosian equations, as well as to system identification (\cite{MQ2}). We note that for the index range $1<p<\infty$ no matter whether ({\ref{non-separate condition}) is met or not the generalized system $\{B_k\}$ ia a Schauder basis of the  $L^p({\bf \partial D})$-topological closure of the ${\rm span}\{B_k\}$  (\cite{QChenTan}).\\

Now we discuss whether the idea od AFD can be generalized to multiple dimensions. AFD is based on complex analysis of one complex variable. For a Euclidean space ${\bf R}^n$ there exist essentially two formulations to bring in a complex structure, or a Cauchy theory. One is the  several complex variables formulation with the imbedding ${\bf R}^n\subset {\bf C}^n.$ It changes a function of several real variables, $f(x_1,...,x_n),$ to the corresponding one in several complex variables, $f(z_1,...,z_n).$ This makes sense at least when $f$ is a polynomial. In such formulation $f$ is said to be holomorphic if and only if $f$ is holomorphic in each component variable $z_k, k=1,2,...,n.$ The other is the imbedding into a Clifford algebra, with the quaternionic space as a particular case. The formulation extends the domain of a function $f(x_1,...,x_n)$ to a set of one more dimension, $f(x_0,x_1,...,x_n).$ With each of those formulations there exists a Hardy space theory. In the several complex variables case there exist Hardy spaces on tubes (\cite{St}, \cite{SW}); while, for the Clifford algebra formulation, there exists a conjugate harmonic system theory (\cite{St}, \cite{SW}), being alternative to the notion of Clifford Hardy space theory (\cite{LMcQ}, \cite{KQ}). The Clifford-quaternionic formulation is more close to one complex variable: 1-D AFD has been extended or partially extended to the quaternionic and the Clifford algebra contexts. The first paper along this line is \cite{QSW} generalizing AFD to quaternions. The second paper is \cite{QWY} dealing with general Clifford algebras. Due to the non-commutativity obstacle what one can do in ${\bf R}^n$ is a greedy-like algorithm. The ADF generalizations to those contexts are not in full due to the non-commutative obstacle. The third paper is \cite{YQS} in which  a scalar-valued phase and its derivative are introduced, and used to analyze signals of several real variables. The phase derivative concept has a close connection with signal decomposition. The success of the generalizations to the quaternionic-Clifford algebra formulation lays on the fact that one can perform algebraic computations to higher dimensional vectors just like what one does to complex numbers. The inconvenience of such formulation is that, except the two cases, viz., the whole hyper-planes and the real spheres in ${\bf R}^n,$ there exist very few cases in which we could possibly use the quaternionic and Clifford algebra theory. For images defined in bounded regions, including rectangular regions for instance, one runs into difficulty. \\

 For the several complex variables formulation there is no direct and parallel method like what we use in the 1-D case to produce an AFD theory. In particular, there is no backward shift mechanism, nor outer functions theory, nor Blaschke products either. The present study offers two approaches for the several complex variables formulation of which one is phrased as $\lq\lq$Product-TM Systems"; and the other $\lq\lq$Product-Szeg\"o Dictionary", as provided, respectively, in \S 2 and \S 3.  With the Product-TM System approach, like in the 1-D theory, we prove a Maximal Selection Principle and then prove the corresponding convergence. With the Product-Szeg\"o Dictionary approach, we first show that Orthogonal Greed Algorithm (OGA) is applicable, ie. the maximal energy may be gained at each step. We next propose a new type greedy algorithm in the abstract complex Hilbert space setting, called Pre-Orthogonal Greedy Algorithm, abbreviated as P-OGA. Our proofs of the corresponding convergence and convergence rate estimation are made to be more general by allowing the weak type version. To be consistent with the literature, by $\lq\lq$weak" we mean to allow a tolerance constant $\rho <1$ in the Maximal Selection Principle (\ref{New}). The new convergence rate estimation (Theorem \ref{rate}) evidences that P-OGA is stronger than OGA. In \S 4 we further study P-OGA in depth and raise a new concept \emph{Induced Complete Dictionary}, or Complete Dictionary.  We show that the concept is naturally associated with P-OGA. We then prove the availability of P-OGA in the 2-torus context under the Complete Dictionary induced by the Product-Szeg\"o Dictionary. Finally we prove that in the classical 1-D unit disc context, under the Complete Dictionary induced from the Szeg\"o Dictionary, P-OGA is identical with AFD.\\

\section {2-D AFD of the Product-TM System Type}

  We will be working with two complex variables. For more several complex variables there is a parallel theory.

Let ${\bf a}$ denote a finite or infinite sequence $\{a_n\}$ of  complex numbers $a_1,a_2,...$ in the unit disc ${\bf D},$ and ${\cal B}^{\bf a}$ the finite or infinite TM system defined through ${\bf a}$ , ie.
\[{\cal B}^{\bf a}=\{B_{\{a_1,...,a_n\}}\}=\{B^{\bf a}_n\},\]
where
\[ B^{\bf a}_n(z)=\frac{\sqrt{1-|a_n|^2}}{1-\overline{a}_nz}\prod_{l=1}^{n-1}\frac{z-a_l}{1-\overline{a}_lz}, \quad n=1,2,...\]
When ${\bf a}$ is a finite sequence, ${\bf a}=\{a_1,...,a_N\},$ we sometimes denote ${\cal B}^{\bf a}$ by
${\cal B}^{\bf a}_N.$

Let ${\bf T}$ denote the boundary of the unit disc $\partial {\bf D},$ and $L^2({\bf T}^2)$ the space of complex-valued functions on the 2-torus with finite energy, where the energy is defined via the inner product
\[\langle f,g\rangle = \frac{1}{4\pi^2}\int_{-\pi}^{\pi}\int_{-\pi}^{\pi} f(e^{it},e^{is})\overline{g}(e^{it},e^{is})dtds.\]

From the multiple trigonometric series theory $f\in L^2 ({\bf T}^2)$ if and only if
\[  f(e^{it}, e^{is})=\sum_{-\infty <k,l< \infty}c_{kl}e^{i(kt+ls)} \qquad {\rm  in\ the\ } L^2-{\rm sense},\]
where
\[ \sum_{-\infty <k,l< \infty} |c_{kl}|^2 < \infty, \quad c_{kl}=\langle f, e_{kl}\rangle,\quad  e_{kl}(t,s)=e^{ikt}e^{ils}.\]

Denote
\[ H^2 ({\bf T}^2)=\{ f\in L^2 ({\bf T}^2) \ : \ f(e^{it}, e^{is})=\sum_{k,l\ge 0}c_{kl}e^{i(kt+ls)}\}.\]
In the sequel if there is no confusion arising we will denote $L^2 ({\bf T}^2)$ and $H^2 ({\bf T}^2),$ briefly and respectively, by $H^2$ and $L^2.$
 It can be easily shown that $H^2$ is a closed subspace of $L^2.$
Denote by $H^2({\bf D}^2)$ the class of complex holomorphic functions in the poly-disc
${\bf D}\times {\bf D}$ satisfying
\[ \sup_{0<r, s<1} \int_{-\pi}^\pi \int_{-\pi}^\pi |f(re^{it}, se^{iu})|^2dtdu < \infty.\]
It may be shown that for any function  $f\in H^2({\bf D}^2)$
there holds
\[ \lim_{z\to e^{it}; w\to e^{is}} f(z,w) \qquad {\rm exist\ for \ almost \ all } \ (e^{it}, e^{is})\in {\bf T}^2,\]
where both the limits $z\to e^{it}$ and $ w\to e^{is}$ are in the non-tangential manner in their respective unit disc, and the limit function on ${\bf T}^2$ belongs to $H^2({\bf T}^2).$ The mapping that maps a function in $f\in H^2({\bf D}^2)$ to its boundary limit function in $H^2({\bf T}^2)$ is one to one and onto, and, as a matter of fact, an isometric isomorphism. For this reason we sometimes use $H^2$ for both  $H^2({\bf D}^2)$ and $H^2({\bf T}^2).$ We note that ${\bf T}^2$ is not the topological boundary, but \emph{characteristic boundary} of ${\bf D}^2,$ which is part of the topological boundary. Data on the characteristic boundary, however, determine a holomorphic function inside ${\bf D}^2$ through the Cauchy integral with the tensor type Cauchy kernel of two complex variables.

 From now on we work with a real-valued function $f\in L^2.$ Define
 \[f^{+,+}(e^{it},e^{is})=\sum_{k,l\ge 0}c_{lk}e^{i(kt+ls)},\]
\[f^{+,-}(e^{it},e^{is})=\sum_{k,-l\ge 0}c_{lk}e^{i(kt+ls)},\]
\[f^{-,+}(e^{it},e^{is})=\sum_{-k,l\ge 0}c_{lk}e^{i(kt+ls)},\]
\[f^{-,-}(e^{it},e^{is})=\sum_{-k,-l\ge 0}c_{lk}e^{i(kt+ls)}.\]

Analogous  with (\ref{real valued 1}) we have

\begin{theorem} Let $f\in L^2$ be real-valued. Then
\[f(e^{it},e^{is})=2{\rm Re}\{f^{+,+}\}(e^{it},e^{is})+2{\rm Re}[f(e^{i( \cdot)}, e^{-i(\cdot)})]^{+,+} (e^{it},e^{-is})-2{\rm Re}\{F^+\}(e^{it})-2{\rm Re}\{G^+\}(e^{is})+c_{00}.\]
\end{theorem}
\noindent{\bf Proof}
We have the relation
\begin{eqnarray*} f(e^{it},e^{is})&+&F(e^{it})+G(e^{is})+c_{00}=f^{+,+}(e^{it},e^{is})+ f^{+,-}(e^{it},e^{is})+\\
& & \qquad +f^{-,+}(e^{it},e^{is})+f^{-,-}(e^{it},e^{is}),\end{eqnarray*}
where
\[ F(e^{it})=\frac{1}{2\pi}\int_{-\pi}^\pi f(e^{it},e^{is})ds, \qquad G(e^{is})=\frac{1}{2\pi}\int_{-\pi}^\pi f(e^{it},e^{is})dt.\]
Therefore,
\begin{eqnarray*} f(e^{it},e^{is})&=& f^{+,+}(e^{it},e^{is})+ f^{+,-}(e^{it},e^{is})+f^{-,+}(e^{it},e^{is})+\\
& & \qquad + f^{-,-}(e^{it},e^{is})-F(e^{it})-G(e^{is})-c_{00}.\end{eqnarray*}
We note that
\[[f(e^{i(\pm \cdot)},e^{i(\pm \cdot)})]^{+,+} (e^{i (\pm t)},e^{i (\pm s)})=f^{\pm ,\pm}(e^{it},e^{is}), \quad  [f(e^{i(\pm \cdot)},e^{i(\mp \cdot)})]^{+,+} (e^{i(\pm t)},e^{i(\mp s)})=f^{\pm,\mp }(e^{it},e^{is}).\]
Since $f$ is real-valued, $f^{+,+}$ and $f^{-,-}$ is a conjugate pair, and  $f^{+,-}$ and $f^{-,+}$ is another conjugate pair.  We have
\[ f^{+,+}+f^{-,-}=2{\rm Re}\{f^{+,+}\}\]
and
\[  f^{+,-}+f^{-,+}=2{\rm Re}\{f^{+,-}\}.\]
Thus,
\begin{eqnarray*} f(e^{it},e^{is})&=& 2{\rm Re}\{f^{+,+}\}(e^{it},e^{is})+2{\rm Re}\{f^{+,-}\}(e^{it},e^{is})-F(e^{it})-G(e^{is})-c_{00}\\
&=&2{\rm Re}\{f^{+,+}\}(e^{it},e^{is})+2{\rm Re}[f(e^{i( \cdot)}, e^{-i(\cdot)})]^{+,+} (e^{it},e^{-is})-F(e^{it})-G(e^{is})-c_{00}\\
&=&2{\rm Re}\{f^{+,+}\}(e^{it},e^{is})+2{\rm Re}[f(e^{i( \cdot)}, e^{-i(\cdot)})]^{+,+} (e^{it},e^{-is})-2{\rm Re}\{F^+\}(e^{it})-2{\rm Re}\{G^+\}(e^{is})+c_{00}.\end{eqnarray*}
The proof is complete.\\

We note that $f^{+,+}(e^{it}, e^{is}),f^{+,-}(e^{it},e^{-is}), f^{-,+}(e^{-it},e^{is}), f^{-,-}(e^{-it},e^{-is})$ are functions in $H^2.$  The above result shows that decomposition of a real-valued function $f\in L^2$ may be
 reduced to decomposition of a number of related functions in the Hardy space.

It may be easily shown that the Product-TM System is complete in the product space, if the two factor 1-D TM systems both are complete in their
 respective spaces.

\begin{theorem}If  ${\cal B}^{\bf a}_N$ and ${\cal B}^{\bf b}_M$ are two finite TM systems, then
${\cal B}^{\bf a}_N \bigotimes {\cal B}^{\bf b}_M$ is an orthonormal system in $L^2({\bf T}^2).$ When
 ${\cal B}^{\bf a}$ and ${\cal B}^{\bf b}$ are two bases of $H^2({\bf T}),$ then ${\cal B}^{\bf a} \bigotimes {\cal B}^{\bf b}$ is a basis of $H^2({\bf T}^2).$
\end{theorem}

\noindent{\bf Proof} The first assertion is obvious. To show the second, we note that holomorphic functions in two complex variables of the type $\sum_{k=1}^K f_k(z)g_k(w)$ is dense in $H^2({\bf T}^2).$  In fact, finite sums of multiple trigonometric series are dense in  $H^2({\bf T}^2).$ If ${\cal B}^{\bf a}$ and ${\cal B}^{\bf b}$ are two bases of $H^2({\bf T}),$ then finite linear combinations of functions in  ${\cal B}^{\bf a} \bigotimes {\cal B}^{\bf b}$ are dense in the function class consisting of functions of the type $\sum_{k=1}^K f_k(z)g_k(w),$ and therefore also   dense in $H^2({\bf T}^2).$  The proof is complete.\\

Denote, for $f\in H^2 ({\bf T}^2),$
\begin{eqnarray}\label{expansion}  S_n(f)=\sum_{1\leq k,l\leq n} \langle f,B^{\bf a}_k\otimes B^{\bf b}_l\rangle B^{\bf a}_k\otimes B^{\bf b}_l&=&\sum_{k=1}^n D_n(f), \ D_n(f)=S_n(f)-S_{n-1}(f), \ S_0(f)=0,\\
& & \qquad \qquad \qquad \qquad \qquad \qquad \qquad n=1,2,... \nonumber
\end{eqnarray}
Note that $D_n(f)$ is called the $n$-\emph{partial sum difference} having $2n-1$ entries.

\begin{theorem}{\rm (Maximal Selection Principle for Product-TM System)} \label{maximal} For  any $f\in H^2$ and
previously fixed $a_1,...,a_{n-1}$ and
$b_1,...,b_{n-1}$ in ${\bf D}$ there exist $a_n,b_n$ in ${\bf D}$ such that the associated
\begin{eqnarray}\label{does not} \| D_n (f) \|^2 =   \sum_{\max\{k,l\}=n} |\langle f,B^{\bf a}_k\otimes B^{\bf b}_l\rangle|^2\end{eqnarray}
attains its maximal value among all possible selections of $a_n, b_n$ inside the unit disc.
\end{theorem}
\noindent{\bf Proof} Let $f\in H^2$ be given and fixed.
We separate the proof into two steps: (i) when $|a_n|\to 1$ and $ |b_n|\to 1,$ with  $a_1,...,a_{n-1}$ and $b_1,...,b_{n-1}$ being previously fixed, one has, uniformly in $a_1,...,a_{n-1}$ and $b_1,...,b_{n-1},$
\[ \lim_{|a_n|\to 1, |b_n|\to 1} \| D_n(f)\|^2=0;\]
and, (ii) if one of $|a_n|$ and $|b_n|$ tends to $1,$ then (\ref{does not}) does not give rise to a maximal value either. We note that the first assertion deals with the case where $(a_n,b_n)$ tends to the characteristic boundary; while the second assertion deals with the boundary.

Now we show (i).
Applying the Cauchy-Schwarz inequality to each of the $2n-1$ terms of the partial sum difference $D_n,$ we conclude that, for any $\epsilon >0,$ we can find a polynomial
$P$ such that
\[ \| D_n (f-P) \|^2 \leq \epsilon\]
uniformly in $a_1,...,a_{n-1}, a_n$ and $b_1,...,b_{n-1}, b_n.$  It is therefore reduced to showing that for any polynomial $P,$
 \[ \lim_{|a_n|\to 1, |b_n|\to 1} \| D_n(P)\|^2=0.\]

 Now due to the expansion of $D_n$ in (\ref{expansion}), it suffices to show, under the limit procedure $|a_n|\to 1 $ and $|b_n|\to 1,$

 \begin{eqnarray}\label{first}
 |\langle P, B^{\bf a}_n\otimes B^{\bf b}_l\rangle|^2 \to 0, \qquad 1\leq l\leq n\end{eqnarray}
and
 \begin{eqnarray}\label{second}
 |\langle P, B^{\bf a}_k\otimes B^{\bf b}_n\rangle|^2 \to 0, \qquad 1\leq k<n.\end{eqnarray}

 Denote by $S_a^{(1)}$ the generalized partial backward shift operator via $a$ for the variable $z,$ and similarly $S_b^{(2)}$ for the second variable $w.$ Due to (\ref{relations}), as well as the reproducing property of the product Szeg\"o kernel $e_{a_n}\otimes e_{b_n},$  we have

 \begin{eqnarray}\label{only}\langle P, B^{\bf a}_n\otimes B^{\bf b}_l\rangle &=& \langle \prod_{k=1}^{n-1}S_{a_k}^{(1)} \prod_{k=1}^{l-1}S_{b_k}^{(2)} (P),
 e_{a_n}\otimes e_{b_l}  \rangle \nonumber \\
 &=& \sqrt{1-|a_n|^2}\sqrt{1-|b_l|^2}\prod_{k=1}^{n-1}S_{a_k}^{(1)} \prod_{k=1}^{l-1}S_{b_k}^{(2)} (P)(a_n,b_l)\nonumber \\
 &\to& 0,\qquad {\rm as} \quad |a_n| \to 1,
 \end{eqnarray}
 where we use the fact that generalized backwards shifts of polynomials are still polynomials, and thus are bounded in a neighborhood of the closed unit disc.
 Similarly we can show (\ref{second}).
 We thus have
 \[  \lim_{|a_n|\to 1, |b_n|\to 1} \| D_n (P) \|^2=0.\]
Now we show (ii). Let $|a_n|\to 1. $ Due to (\ref{only}), all terms in (\ref{first}) are vanishing. The orthogonality between the $2n-1$ terms of $D_n(f)$ implies that for $|a_n|\to 0$ only the terms corresponding to (\ref{second}) can contribute to the energy gain, and thus do not give rise to the maximal energy, unless for all $a_n$ the terms in (\ref{1st}) are zero. In the latter case the remainder after the former $n-1$ steps depends only on the variable $w,$ reducing to one complex variable case. In any case the maximal energy is attained only at an interior point.
The proof is complete.\\

Thanks to orthogonality between distinct $B^{\bf a}_k\otimes B^{\bf b}_l$'s, there follows
\[ 0\leq \| f-\sum_{k=1}^n D_k(f) \|^2 = \|f \|^2 - \sum_{k=1}^n \| D_k(f)\|^2\]
that implies the Bessel inequality
\[ \sum_{k=1}^n \| D_k(f)\|^2\leq \|f\|^2.\]
As consequence, we have
\[ \lim_{n\to \infty} \sum_{k=n+1}^\infty \| D_k(f)\|^2 =0.\]
This, in particular, is valid with maximal selections of $(a_n, b_n)$ in accordance with Theorem \ref{maximal}.
Moreover, under such selections we have

\begin{theorem}\label{Th3} Let $f\in H^2 ({\bf T}^2).$ For any $k_0$ and previously fixed $a_1,b_1,...,a_{k_0-1}, b_{k_0-1},$ by selecting parameter pairs  $(a_{k_0}, b_{k_0}), (a_{k_0+1}, b_{k_0+1}), ..., $ according to the Maximal Selection Principle (MSP), we have
\[  \lim_{n\to \infty} \| f - S_n(f) \|^2 = 0. \]
 In other words, in the $L^2$-convergence  sense,
\[  f=\lim_{n\to \infty} S_n(f). \]
\end{theorem}
\noindent{\bf Proof} We prove this by contradiction. Assume that this is not true. Then
\[ f=\sum_{k=1}^\infty D_k(f)+h, \quad h\ne 0,\]
where $h$ is in $H^2,$ and orthogonal with each $D_k(f).$ Hence,
\[ \|h\|^2=\| f\|^2 - \sum_{k=1}^\infty \| D_k(f)\|^2 > 0. \]
 By using the tensor type Cauchy formula for two complex variables we have for any  $\tilde{a}, \tilde{b}$ in ${\bf D},$
 \[  \langle h, e_{\{\tilde{a}\}}\otimes e_{\{\tilde{b}\}}\rangle = \sqrt{1-|\tilde{a} |^2}\sqrt{1-|\tilde{b} |^2}h(\tilde{a}, \tilde{b})\]
Therefore, there exist $\tilde{a}, \tilde{b}$ in the unit disc such that $\langle h, e_{\{\tilde{a}\}}\otimes e_{\{\tilde{b}\}}\rangle \ne 0.$
Denote by
\[ \tilde{X}=\overline{{\rm span}}\{\tilde{\cal B}^{\bf a}\otimes\tilde{\cal B}^{\bf b}\},\]
where $\tilde{\cal B}^{\bf a}$ is the TM system generalized by $\{ \tilde{a}, a_1, ... , a_n, ...\}$ under the given order; and likewise for the notation $\tilde{\cal B}^{\bf b}.$ The corresponding $n$-partial sum difference is denoted
$\tilde{D}_n,$ involving $\tilde{a}, \tilde{b}, a_1, b_1,..., a_{n-1}, b_{n-1}.$ \\

Denote by $h/\tilde{X}$ the orthogonal projection of $h$ into the subspace $\tilde{X}.$ In the sequel we will continue to adopt this notation for orthogonal projections. It is easy to show that
\[ \| h/\tilde{X} \|_2=\delta >0.\]
In fact,
\[ \| h/\tilde{X} \|_2^2\ge \sum_{k=1}^\infty \| \tilde{D}_k \|_2^2 \ge \| \tilde{D}_1 \|_2^2 = | \langle h, e_{\tilde{a}}\otimes e_{\tilde{b}}\rangle |^2 >0.\]
Set, for any integer $M,$
\[ \tilde{X}_M={\rm span} \{ \tilde{\cal B}^{\bf a}_M \otimes \tilde{\cal B}^{\bf b}_M\}\quad {\rm and}\quad {X}_M={\rm span} \{ {\cal B}^{\bf a}_M \otimes {\cal B}^{\bf b}_M\},\]
where $\tilde{\cal B}^{\bf a}_{M}$ and  $ \tilde{\cal B}^{\bf b}_{M}$ are, respectively, the
TM systems generalized by $\{\tilde{a}, a_1, ... , a_{M-1}\}$ and $\{\tilde{b}, b_1, ... , b_{M-1}\}.$

Since
\[ \| h/\tilde{X} - h/\tilde{X}_M \|^2=\sum_{k=M+2}^\infty \| \tilde{D}_k \|^2 \to 0,\]
we have
\[ \lim_{M\to \infty} h/\tilde{X}_M = h/ \tilde{X}.\]
Now, find $M$ so large that
\[\| h/\tilde{X}_M \|^2 > \delta /2 \qquad {\rm and} \qquad \| \sum_{k=M}^\infty D_k(f) \|^2 < \delta /8.\]
For such $M,$ on one hand,
\[ \| f/X_{M} \|^2 = \| \sum_{k=1}^{M-1} D_k(f) \|^2 + \| D_{M}(f) \|^2 < \| \sum_{k=1}^{M-1} D_k(f) \|^2 + \delta / 8,\]
that is for a maximal selection of $(a_{M}, b_{M}).$
On the other hand,
 \[ \| f/\tilde{X}_M \|^2 = \| \sum_{k=1}^{M-1} D_k(f) \|^2 + \| \left(\sum_{k=M}^\infty D_k(f)\right)/ \tilde{X}_M + h/\tilde{X}_M  \|^2>\| \sum_{k=1}^{M-1} D_k(f) \|^2+\delta /2 - \delta /8.\]
 This shows that the selection $(\tilde{a}, \tilde{b})$ is better than the selection $(a_{M}, b_{M}),$ being contradictory to the maximality of the pair $(a_{M},b_{M}).$
 The proof of Theorem \ref{Th3} is, therefore, complete.\\

The approximation algorithm in Theorem \ref{Th3} is called \emph{Two Dimensional AFD} (\emph{2D-AFD}) \emph{of the Product-TM System type}.\\

\section {2-D AFDs of the Product-Szeg\"o Dictionary type}

We, in general, call all types of adaptive approximation by linear combinations of parameterized Szeg\"o and higher order Szeg\"o kernels AFDs.  We denote by
\[ {\cal D}^2=\{ e_{a}\otimes e_{b} \ : \ a,b\in {\bf D}\}\] the set of functions consisting of tensor products of two 1-D parameterized and normalized Szeg\"o kernels.
It is easy to show that ${\cal D}^2$ is a dictionary of the $H^2$ space of the poly-disc.
We call it \emph{Product-Szeg\"o Dictionary}. We first verify that with the dictionary ${\cal D}^2$ the Pure Greedy Algorithm (PGA) (\cite{Tem}) is applicable to $H^2.$
Write
\[ f(z,w)=\sum_{k=1}^{n-1} \langle g_k,  e_{a_k}\otimes e_{b_k}\rangle e_{a_k}(z)e_{b_k}(w)+\tilde{g}_{n}(z,w), \qquad \tilde{g}_1=f.\]
We note that the \emph{standard remainders} $\tilde{g}_{n}$ defined above is different from what is defined in (\ref{standard}), the latter, being phrased as \emph{orthogonal standard remainders}, related to the orthogonal projection of the given function to the linear space generated by the selected dictionary elements.

Since $\tilde{g}_{k+1}$ is orthogonal with $e_{a_k}\otimes e_{b_k},$ we have the relation
\begin{eqnarray}\label{repeaking} \| \tilde{g}_{k+1} \|^2 = \| \tilde{g}_k-\langle \tilde{g}_k, e_{a_k}\otimes e_{b_k}\rangle e_{a_k}\otimes e_{b_k} \|^2=\|\tilde{g}_k\|^2-| \langle \tilde{g}_k, e_{a_k}\otimes e_{b_k} \rangle |^2.\end{eqnarray}
The associated energy rule is
  \[ \| f- \sum_{k=1}^n \langle \tilde{g}_k, e_{a_k}\otimes e_{b_k} \rangle (e_{a_k}\otimes e_{b_k})\|^2=\|f\|^2-\sum_{k=1}^n
  |\langle \tilde{g}_k, e_{a_k}\otimes e_{b_k} \rangle |^2\]
that is obtained through recursive application of (\ref{repeaking}).

We show that for each $k$ the pair $(a_k, b_k)$ can be selected such that
\[ (a_k,b_k)=\arg \max \{ |\langle \tilde{g}_k, e_{a}\otimes e_{b}\rangle| \ : \ a,b\in {\bf D}\}.\]

\begin{theorem}\label{Maximal Szego}
{\rm (Maximal Selection Principle for Product-Szeg\"o Dictionary)} \label{maximal2} For  any $\tilde{g}\in H^2({\bf D}^2)$ one can find
\[ (\tilde{a},\tilde{b})=\arg \max \{ |\langle \tilde{g}, e_{a}\otimes e_{b}\rangle| \ : \ a,b\in {\bf D}\}.\]
\end{theorem}

The proof of Theorem \ref{maximal} by using polynomial approximation
can be easily adapted to give a proof of Theorem \ref{maximal2}. We, however, give an alternative proof that, by the author's opinion, has advantages in proving existence of the global maximum under simultaneous selections of several parameters (see \cite{MQW}). The idea of the following proof
 was originated by Temlyakov (see \cite{Te}, \cite{QWa}).\\

\noindent {\bf Proof} It suffices to show for any function $\tilde{g}$ in the Hardy space there holds
\[ \lim_{|a|\to 1- \ {\rm or}\  |b|\to 1-} |\langle \tilde{g}, e_{a}\otimes e_{b}\rangle |^2=0.\]
Due to the orthogonality it suffices to show
\begin{eqnarray}\label{lim} \lim_{|a|\to 1- \ {\rm or}\  |b|\to 1-} \| \tilde{g}-\langle \tilde{g}, e_{a}\otimes e_{b}\rangle e_{a}\otimes e_{b} \| =\| \tilde{g} \|.\end{eqnarray}
Let $P_r\otimes P_s$ be the tensor type Poisson kernel on the poly-disc (\cite{SW}), $r,s\in [0,1).$ For $\epsilon >0,$ we can choose $r$ and $s$ sufficiently close to $1$ such that, owing to the $L^2$-approximation property of the Poisson integral, we have
\begin{eqnarray}\label{inequality}
\| \tilde{g}\| &\ge& \| \tilde{g}-\langle \tilde{g}, e_{a}\otimes e_{b}\rangle e_{a}\otimes e_{b} \| \nonumber\\
& \ge & \| (P_r\otimes P_s) \ast [\tilde{g}-\langle \tilde{g}, e_{a}\otimes e_{b}\rangle e_{a}\otimes e_{b}]\|\nonumber \\
& \ge & \|  (P_r\otimes P_s) \ast \tilde{g} \| - |\langle \tilde{g}, e_{a}\otimes e_{b}\rangle  | \| (P_r\otimes P_s) \ast (e_{a}\otimes e_{b})\|\nonumber \\
& \ge & (1-\epsilon )\|\tilde{g}\| - \| \tilde{g}\| \|  (P_r\otimes P_s) \ast (e_{a}\otimes e_{b})\|.\end{eqnarray}
Now, with the fixed $r$ and $s$, since $e_{a}\otimes e_{b}\in H^2,$ there follows, for $z=re^{it}, w=se^{iu},$
\[ (P_r\otimes P_s) \ast (e_{a}\otimes e_{b})(e^{it},e^{iu})=e_{a}(z)e_{b}(w).\]
Then we have explicit computation
\begin{eqnarray*}
\| (P_r\otimes P_s) \ast (e_{a}\otimes e_{b})\|^2 &=& \frac{1}{(2\pi )^2} \int_0^{2\pi} \frac{1-|a|^2}{|1-\overline{a}r e^{it}|^2}dt\int_0^{2\pi} \frac{1-|b|^2}{|1-\overline{b}s e^{iu}|^2}du\\
&=& \frac{1-|a|^2}{1-r^2|a|^2}\frac{1-|b|^2}{1-s^2|b|^2}.\end{eqnarray*}
When $|a|\to 1$ or $ |b|\to 1,$ the inequality (\ref{inequality}) gives
\[ \|\tilde{g}\|\ge \| \tilde{g}-\langle \tilde{g}, e_{a}\otimes e_{b}\rangle e_{a}\otimes e_{b} \|\ge (1-2\epsilon )\|\tilde{g}\|.\]
This shows that the limit (\ref{lim}) holds. The proof is complete.

 By implementing the above Maximal Selection Principle with the Product-Szeg\"o Dictionary, under the general theory of PGA (\cite{MaZh}, \cite{DT}, \cite{DMA}, \cite{Tem}),
  one has, in the energy sense,
  \[ f=\sum_{k=1}^\infty \langle g_k, e_{a_k}\otimes e_{b_k} \rangle e_{a_k}\otimes e_{b_k}.\]\\

  It was shown that Orthogonal Greedy Algorithm (OGA) is more optimal than PGA (\cite{DT}, \cite{DMA}, \cite{Tem}). Below we propose a new variation of greedy algorithm called Pre-Orthogonal Greedy Algorithm (P-OGA) that is truly stronger than OGA. Our formulation is in the general complex Hilbert space, allowing a weak type version of P-OGA, called WP-OGA.

Let ${\cal H}$ be a complex Hilbert space and ${\cal A}$ a dictionary consisting of elements  $a\in {\cal A}$ satisfying $\| e_a\|=1, \ \overline{{\rm span}\cal A}={\cal H}.$ Let
\begin{eqnarray}\label{let} f=\sum_{k=1}^{n-1} \langle f, B_k\rangle B_k + g_{n},\end{eqnarray}
where $\{ B_1,...,B_k\}$ is the result of the Gram-Schmidt orthogonalization
 process applied to the finite system $\{ a_1,...,a_k \},$ where for each $k,$
$a_k$ is selected to be one satisfying the \emph{Pre-Orthogonal $\rho$-Maximal Selection Principle}
\begin{eqnarray}\label{New} |\langle g_k, B_k \rangle | \ge \rho \sup \{ |\langle g_k, B_k^a \rangle | \ :\ a\in {\cal A}\}, \quad \rho\in (0,1],\end{eqnarray}
 where $\{B_1,...,B_{k-1},B_k^a\}$ is the orthogonalization of $\{a_1,...,a_{k-1},a\}.$  Such defined decomposition is called \emph{Weak Pre-Orthogonal Greedy Algorithm}, abbreviated as WP-OGA. When $\rho=1,$ it is called \emph{Pre-Orthogonal Greedy Algorithm}, abbreviated as \emph{P-OGA}. To see the difference between the just introduced with \emph{Weak Orthogonal Greedy Algorithm } (\emph{WOGA}, see \cite{Tem}) we recall that in the latter the selected $a_k$ satisfies
\begin{eqnarray}\label{Old} |\langle g_k, a_k \rangle | \ge \rho \sup \{ |\langle g_k, a \rangle | \ :\ a\in {\cal A}\}, \quad \rho\in (0,1].\end{eqnarray}
 To summarize, both WP-OGA (P-OGA) and WOGA (OGA) use the orthogonal standard remainder but WP-OGA first performs orthogonalization and then selects a dictionary element, while WOGA first selects a dictionary element and then performs orthogonalization. To show that WP-OGA offers a better approximation than OGA at each step, we note that, with a selected $a_n\in {\cal A},$ no matter how  $| \langle g_k, a_n \rangle |$ is close to $\sup \{ |\langle g_k, a \rangle | \ :\ a\in {\cal A}\},$ we always have
 \[ \sup \{ |\langle g_k, B_k^a \rangle | \ :\ a\in {\cal A}\} \ge |\langle g_k, B^{a_n}_n \rangle|.\]
The optimality of WP-OGA is also seen from the convergence rate estimation proved in Theorem \ref{rate}.\\

As in OGA, we have, owing to the orthogonality,
  \[ \langle f, B_k \rangle = \langle g_k, B_k \rangle.\] The corresponding energy rule is again
 \[ \| f- \sum_{k=1}^n \langle f, B_k\rangle B_k\|^2=\|f\|^2-\sum_{k=1}^n
  |\langle f, B_k \rangle |^2\]
that implies the Bessel type inequality
  \[ \sum_{k=1}^\infty
  |\langle f, B_k \rangle |^2\leq \| f\|^2.\]

\begin{theorem}\label{prove} Let ${\cal H}$ be a complex Hilbert space with a dictionary ${\cal A}.$  For any $f\in {\cal H},$
 with a sequence of consecutively selected $a_1,...,a_n,...$ from ${\cal A}$ under the Pre-Orthogonal $\rho$-Maximal Selection Principle  we have
\[ f=\sum_{k=1}^\infty \langle f, B_k\rangle B_k,\]
where for each $k$ the system $\{B_1,...,B_k\}$ is the result of the Gram-Schmidt orthonormalization process applied to $\{a_1,...,a_k\}.$
\end{theorem}
The proof provided below directly depends on the Pre-Orthogonal $\rho$-Maximal Selection Principle.  \\

\noindent{\bf Proof} We prove the theorem by contradiction. Assume that $f=\sum_{k=1}^\infty \langle f, B_k\rangle B_k+h,$
$h\ne 0, h\perp \overline{{\rm span}{\cal B}},$ where ${\cal B}={\rm span}\{a_1,...,a_n,...\}.$ Since $\overline{{\rm span}\cal A}={\cal H},$ there exists $b\in {\cal A}$ such that $\langle h, e_b\rangle \ne 0.$ Denote ${\cal B}^b=\{ b, a_1,...,a_n,...\}.$ Obviously, $h/\overline{{\rm span}{\cal B}^b}\ne 0.$ Let $\| h/\overline{{\rm span}{\cal B}^b}\|=\delta \ (>0).$ Denote ${\cal B}_n=
\{a_1,...,a_n\}, {\cal B}^b_{n+1}=\{ b,a_1,...,a_n\}.$ We have, by similar reasoning as before,
\[ \lim_{n\to \infty} h/{{\rm span}{\cal B}_{n+1}^b}=h/\overline{{\rm span}{\cal B}^b}.\] Fix $N$ large enough so that
\[ \| h/{{\rm span}{\cal B}_{n+1}^b}\| > \delta /2, \quad {\rm and}\quad \| \sum_{k=N+1}^\infty \langle f, B_k\rangle B_k \|< \delta/2^m,\]
where $m$ will be determined later. Now write
\begin{eqnarray*}
f&=&\sum_{k=1}^N \cdot \ + \sum_{k=N+1}^\infty \cdot \ + h\\
&=&\sum_{k=1}^N \cdot \ + g_{N+1}.\end{eqnarray*}
On one hand, due to the orthogonality,
\begin{eqnarray*}
|\langle f, B_{N+1} \rangle | &=& |\langle g_{N+1}, B_{N+1} \rangle |\\
&=& |\langle \sum_{k=N+1}^\infty \cdot \  , B_{N+1} \rangle |\\
&\leq & \| \sum_{k=N+1}^\infty \cdot \ \|\\
&\leq & \delta/2^m.\end{eqnarray*}
On the other hand, since $\{B_1,...,B_N,B^b_{N+1}\}$ is the orthogonalization of $\{ {a_1},...,{a_N}, b\},$ we have
\begin{eqnarray*}
|\langle f, B^b_{N+1}\rangle| &=& |\langle g_{N+1}, B^b_{N+1} \rangle |\\
&=& |\langle h+\sum_{N+1}^\infty \cdot \ , B^b_{N+1}\rangle |\\
&\ge& |\langle h, B^b_{N+1}\rangle |-|\langle \sum_{N+1}^\infty \cdot \ , B^b_{N+1}\rangle |\\
&=& \| h/{{\rm span}{\cal B}_{n+1}^b}\|-|\langle \sum_{N+1}^\infty \cdot \ , B^b_{N+1}\rangle |\\
&\ge& \delta/2 - \delta/2^m\\
&=& \frac{(2^{m-1}-1)\delta}{2^m}.
\end{eqnarray*}
Therefore,
\begin{eqnarray*}
|\langle f, B_{N+1} \rangle |/\sup \{ |\langle f, B_{N+1}^a \rangle | \ :\ a\in {\cal A}\}< \frac{1}{2^{m-1}-1}.\end{eqnarray*}
Now choosing $m$ so large that $ \frac{1}{2^{m-1}-1}<\rho,$ we arrive at a contradiction. The proof is complete.\\

To obtain the convergence rate estimation we first give some remarks. Assume that we have $\{B_1,...,B_{n-1}\}$ as the orthogonalization of $\{{a_1},...,{a_{n-1}}\}.$ When we have the next element ${a_n}$ to be added to, and to be made orthogonal with the former $B_1,...,B_{n-1},$ what we do is to expand ${a_n}$ into the linear expansion $\sum_{k=1}^{n-1} \langle {a_n}, B_k \rangle B_k,$ and then find the $n$-th orthogonal standard remainder ${a_n}- \sum_{k=1}^{n-1} \langle {a_n}, B_k \rangle B_k$ given by (\ref{let}).  This process gives the projection of $a_n$ into the orthogonal complement space of the linear space ${\rm span}\{{a_1},...,{a_{n-1}}\}.$ We denote such projection operator by $Q_{\{{a_1},...,{a_{n-1}}\}}.$ That is,
\[ Q_{\{{a_1},...,{a_{n-1}}\}}({a_n})={a_n}- \sum_{k=1}^{n-1} \langle {a_n}, B_k \rangle B_k, \quad B_n=\frac{Q_{\{{a_1},...,{a_{n-1}}\}}({a_n})}{\| Q_{\{{a_1},...,{a_{n-1}}\}}({a_n})\|}.\]
In below we will sometimes abbreviate $Q_{\{{a_1},...,e_{{n-1}}\}}$ as $Q_{n-1}.$

On the other hand, we also have
\begin{eqnarray}\label{with} g_n=Q_{\{{a_1},...,{a_{n-1}}\}}(f),\end{eqnarray}
where $g_n$ is the $n$-th orthogonal standard remainder of $f$
with respect to the orthonormal system $\{B_1,...,B_{n-1}\}$ defined by (\ref{let}).
The relation we want to cite is that
 for any two functions $f$ and $g$ in the complex Hilbert space we have $\langle f, Q_{\{{a_1},...,{a_{n-1}}\}}(g)\rangle=\langle Q_{\{{a_1},...,e_{{n-1}}\}}(f), g \rangle$ due to the fact that both the left and right hands are identical with
$\langle f, g\rangle - \sum_{k=1}^{n-1} \langle f, B_k \rangle \langle B_k, g \rangle.$

By the just mentioned property of the projection operator, when $\|Q_{\{{a_1},...,{a_{n-1}}\}}(a)\|\ne 0,$ we have
\begin{eqnarray}\label{nonumber}
 |\langle g_n, B_n^a \rangle |&=& \frac{1}{\|Q_{\{{a_1},...,{a_{n-1}}\}}(a)\|} |\langle Q_{\{a_1,...,a_{n-1}\}}(f), Q_{a_1,...,a_{n-1}}(a) \rangle |\nonumber \\
 &=& \frac{1}{\|Q_{\{{a_1},...,{a_{n-1}}\}}(a)\|} |\langle Q^2_{\{a_1,...,a_{n-1}\}}(f), a \rangle |\nonumber \\
 &=&\frac{1}{\|Q_{\{{a_1},...,{a_{n-1}}\}}(a)\|}  |\langle Q_{\{a_1,...,a_{n-1}\}}(f), a \rangle |\nonumber \\
 &=&\frac{1}{\|Q_{\{{a_1},...,{a_{n-1}}\}}(a)\|}  |\langle g_n, a \rangle |.
 \end{eqnarray}
Set, for any $a_1,...,a_{n-1}$ in ${\cal A},$ and any $a\in {\cal A},$
 \begin{eqnarray}\label{as} r_n(a)=\|Q_{\{{a_1},...,{a_{n-1}}\}}(a)\|.\end{eqnarray}
 We have $r_n (a)\leq 1.$ The case $r_n(a)=1$ is exactly when $a$ is orthogonal with ${\rm span}\{{a_1},...,{a_{n-1}}\}.$
 The case $r_n(a)=0$ is exactly when $a\in {\rm span}\{{a_1},...,{a_{n-1}}\}.$ Weak Orthogonal $\rho$-Maximal Selection Principle, however, does not allow such selection. Therefore, for the allowed cases, there always
\[ |\langle g_n, B_n^a \rangle |= \frac{1}{r_n(a)}|\langle g_n, a \rangle | \ge |\langle g_n, a \rangle |.\]

 Now, as in \cite{DT} (also see \cite{QWang}), we introduce
\[ H^2({\cal A}, {M})=\{ f\in H^2 \ : \ f=\sum_{k=1}^\infty c_k a_k, \ \sum_{k=1}^\infty
|c_k|\leq M\}.\]

We have

\begin{theorem}\label{rate} Let $f\in H^2({\cal A}, {M}).$ Denote by $g_m$ the orthogonal standard remainder of $f$ with respect to the orthogonalization $\{B_1,...,B_{m-1}\}$ of the consecutive selected $\{{a_1},...,{a_{m-1}}\}$ under the Pre-Orthogonal $\rho$-Maximal Selection Principle. With $R_m=\max \{r_1,...,r_m\}, \ r_n=\sup_{k} \{ r_n(b_k)\}, $ where $r_n(b_k),$ depending on $a_1,...,a_{n-1}$ and $b_k,$ is defined as in (\ref{as}), we have
\[ \| g_m \|\leq \frac{R_m M}{\rho}\frac{1}{\sqrt{m}}.\]
\end{theorem}

We also need the following result (see \cite{DT})

\begin{lemma}\label{invoking} Let $\{d_n\}_{n=1}^m$ be an $m$-tuple of nonnegative numbers satisfying
\[ d_1\leq A_m, \qquad d_{n+1}\leq d_n\left( 1- \frac{d_n}{A_m}\right).\]
Then there holds
\[ d_m\leq \frac{A_m}{m}.\]
When the above relations hold for all integers $m$ and all $n\leq m,$ and $A_m\leq A,$
then we have, for all $m,$
\[ d_m\leq \frac{A}{m}.\]\end{lemma}

\noindent{\bf Proof of Theorem \ref{rate}} Assume that $f=\sum_k c_k{b_k}$ with $\sum_k |c_k|\leq M.$ We first note that
\[ \| g_{m+1} \|^2=\| g_m \|^2 - |\langle g_m, B_m \rangle |^2.\]
Next, we have a chain of equality and inequality relations: For each $n\leq m,$
\begin{eqnarray*}
| \langle g_n, B_n \rangle | & \ge & \rho \sup_{a\in {\cal A}} | \langle g_n, B_n^a \rangle |\\
&\ge &\rho \sup_{k} | \langle g_n, B_n^{b_k} \rangle |\\
& = &  \rho \sup_{k} \frac{| \langle g_n, b_k \rangle |}{r_n(b_k)}\\
& \ge & \frac{\rho}{r_n} \sup_{k} |\langle g_n, {b_k} \rangle |\\
& \ge & \frac{\rho }{r_n M} |\langle g_n, \sum_k c_k{b_k} \rangle|\\
& = & \frac{\rho}{r_n M} |\langle g_n, f \rangle|\\
& \ge & \frac{\rho }{R_m M} \| g_n \|^2.
\end{eqnarray*}
Therefore,
\[ \| g_{n+1} \|^2\leq \| g_n \|^2 \left(1-(\frac{\rho }{R_m M})^2 \| g_n \|^2\right), \qquad n\leq m.\]
By invoking Lemma \ref{invoking} we obtain the desired estimate. The proof is complete.\\

 Comparing the convergence rate of WOGA (\cite{DT}, \cite{Tem}) with that of WP-OGA given in Theorem \ref{rate}, the latter has an extra factor $Rm.$ For small $m$ the number $R_m$ should be much less than $1.$ Moreover, at each step of selection if the candidates for $a_n$ are multiple, we should select one of those giving rise to the smallest $r_n(a_n), n\leq m,$ and thus the smallest $R_m,$ too.  Under such strategy WP-OGA (P-OGA) is anticipated to achieve considerably better approximation than WOGA (OGA).\\

\section {Complete Dictionary Induced by P-OGA}

In this section we will concentrate in a further study on P-OGA. In an OGA repeated selections of dictionary elements are excluded by its maximal principle. In contrast, a P-OGA with a dictionary ${\cal A}$ whose elements are smoothly parameterized often induces an enlarged dictionary, call \emph{Induced Complete Dictionary}, or \emph{Complete Dictionary} in short, denoted by $\tilde{\cal A}.$ The Maximal Selection Principle of such P-OGA should, in fact, be with respect to the Complete Dictionary $\tilde{\cal A}$ but not ${\cal A}$ (see Theorem \ref{P-OGA applicability} below). The formulation of $\tilde{\cal A}$ is described as follows.

 Now assume that elements of the dictionary ${\cal A}$ are smoothly parameterized by an $m$-dimensional complex vector $T=(t^1,...,t^m)$ in some open set $\Omega \subset {\bf C}^m, m\ge 1.$ In P-OGA, at each step, one encounters the following circumstance: One already has a finite orthonormal system $\{B_1,...,B_{n-1}\}$ being the G-S orthonormalization of the $(n-1)$-tuple $\{a_{T_1},...,a_{T_{n-1}}\}.$ Now with an orthogonal remainder $g_n$ one needs to further select a dictionary element $\tilde{a}=a_{T_n}$ such that the largest possible $|\langle g_n, B^{{a}}_n \rangle|$ is attainable: \begin{eqnarray}\label{assume}|\langle g_n, B^{\tilde{a}}_n \rangle|=\max \{ |\langle g_n, B_n^a \rangle | \ :\ a\in {\cal A}\},\end{eqnarray} where for any $a\in {\cal A}$ the function $B_n^a$ is defined through the requirement that the finite system $\{B_1,...,B_{n-1},B^{{a}}_n\}$ is the G-S orthonormalization of $\{a_{T_1},.,,,a_{T_{n-1}}, a\}.$ If such $\tilde{a}$ exists, then we are done. The point is that such element $\tilde{a}$ may not exist in ${\cal A}.$ We show that under some circumstance it exists in an enlarged dictionary, called the Induced Complete Dictionary, denoted by $\tilde{\cal A}.$ From now on we assume that there does not exist $\tilde{a}\in {\cal A}$ that makes (\ref{assume}) to hold. In the case,
  due to the boundedness $|\langle g_n, B^{\tilde{a}}_n \rangle|\leq \| g_n\|,$ there exists a sequence $T^{(k)}$ converging to a boundary point of $\Omega$ (including $\infty$) or an interior point ${T_n}\in \Omega$ such that
 \begin{eqnarray}\label{interior point}
 \lim_{k\to \infty} |\langle g_n, B^{a_{T^{(k)}}}_n\rangle |\to \sup \{ |\langle g_n, B_n^a \rangle | \ :\ a\in {\cal A}\}|.\end{eqnarray}
 We note that for each $k$ the point $a_{T^{(k)}}$ must not be in the linear span of $a_{T_1},...,a_{T_{n-1}}.$ If it were, then
 \[ Q_n(a_{T^{(k)}})=0,\]
 and hence the new orthonormal function $B^{a_{T^{(k)}}}_n$ cannot be formulated, and there is no further energy gain from
 $|\langle g_n, B^{a_{T^{(k)}}}_n\rangle|.$

 In the classical cases, including the Szeg\"o Dictionary and the Product-Szeg\"o Dictionarie associated with, respectively, the Hardy space in the dist ${\bf D}$ (\cite{QWa}) and one in the poly-disc ${\bf D}^2$ (Theorem \ref{Maximal Szego}),  when $a_{T^{(k)}}$ tends to the boundary of $\Omega$ there holds
 \begin{eqnarray}\label{boundary point}
 \lim_{k\to \infty} |\langle g_n, B^{a_{T^{(k)}}}_n\rangle |\to 0.\end{eqnarray} In the Clifford algebra and the quaternionic cases we observe the same phenomenon (\cite{QSW}, \cite{QWY}). In such cases one concludes that
  the limit point ${T_n}$ of $T^{(k)}$ with the property (\ref{interior point}) must be an interior point of $\Omega .$ We show that $a(T_n)$ must be in ${\rm span}\{a_{T_1},...,a_{T_{n-1}}\}.$ In fact, if not, then we can take $\tilde{a}=a(T_n)$ and (\ref{assume}) holds, being contrary with our assumption. From now on we are based on $a_{T_n}\in {\rm span}\{a_{T_1},...,a_{T_{n-1}}\}.$
 We show that in the case the limit procedure induces a new system function $B^{(\partial_v a_T)({T_n})}_n$ such that
$\{ B_1,...,B_{n-1},B^{(\partial_v a_T)({T_n})}_n\}$ is the orthogonalization of $\{a_{T_1},...,a_{T_{n-1}},
(\partial_v a_T)({T_n})\},$ where $\partial_v$ is the directional derivative along some direction $v, \| v\|=1,$ $\partial_v=v_1\partial_{\overline{t^1}}+\cdots +v_m \partial_{\overline{t^m}}, v=(v_1,...,v_m),
 \partial_{\overline{t^k}}=\frac{1}{2}(\partial_{x_k}-i\partial_{y_k}), t^k=x_k+iy_k, k=1,...,m.$
The assertion is proved through the the limit process
\begin{eqnarray*}& &\lim_{a_{T^{(k)}}\notin {\rm{span}\{a_{T_1},...,a_{T_{n-1}}\}, T^{(k)}\to T_n}}B^{a_{T^{(k)}}}_n \\ &=
  &\lim_{a_{T^{(k)}}\notin {\rm{span}\{a_{T_1},...,a_{T_{n-1}}\}, T^{(k)}\to T_n}}\ \frac{a_{T^{(k)}}-\sum_{k=1}^{n-1}\langle a_{T^{(k)}}, B_k\rangle B_k}{\|a_{T^{(k)}}-\sum_{k=1}^{n-1}\langle a_{T^{(k)}}, B_k\rangle B_k \|}\\
&=&\lim_{a_{T^{(k)}}\notin {\rm{span}\{a_{T_1},...,a_{T_{n-1}}\}, T^{(k)}\to T_n}}
\frac{(a_{T^{(k)}}-\sum_{k=1}^{n-1}\langle a_{T^{(k)}}, B_k\rangle B_k)-({a}_{T_n}-\sum_{k=1}^{n-1}\langle {a}_{T_n}, B_k\rangle B_k)}{\|(a_{T^{(k)}}-\sum_{k=1}^{n-1}\langle a_{T^{(k)}}, B_k\rangle B_k)-({a_{T_n}}-\sum_{k=1}^{n-1}\langle {a_{T_n}}, B_k\rangle B_k) \|}\\
&=&\lim_{a_{T^{(k)}}\notin {\rm{span}\{a_{T_1},...,a_{T_{n-1}}\}, T^{(k)}\to T_n}}
\frac{(a_{T^{(k)}}-{a_{T_n}})/\| T^{(k)}-T_n\|-\sum_{k=1}^{n-1}\langle (a_{T^{(k)}}-{a_{T_n}})/\| T^{(k)}-T_n\|, B_k\rangle B_k)}{\|(a_{T^{(k)}}-{a_{T_n}})/\| T^{(k)}-T_n\|-\sum_{k=1}^{n-1}\langle (a_{T^{(k)}}-{a_{T_n}})/\| T^{(k)}-T_n\|, B_k\rangle B_k)\|}\\
&=& \frac{(\partial_v a_T)(T_n)-\sum_{k=1}^{n-1}\langle (\partial_v a_T) (T_n), B_k\rangle B_k}{\|(\partial_v a_T) (T_n)-\sum_{k=1}^{n-1}\langle (\partial_v a_T (T_n)), B_k\rangle B_k)\|},
\end{eqnarray*}
where $v$ is the tangential direction of the passage $T^{(k)}\to T_n.$ \\

If the parameter space is ${\bf C},$ the above directional derivative reduces to the standard complex derivative, and likewise in the Clifford and quaternionic cases.\\

Denote by ${\cal A}_k, k=1,2,...,$ the function set consisting of all possible normalized directional derivatives of the functions in ${\cal A}_{k-1},$ where ${\cal A}_0={\cal A}.$ Note that for each $k,$ ${\cal A}_k$ is smoothly parameterized by the same parameters $T$ in ${\bf C}^m.$\\

\noindent{\bf Definition}(Induced Complete Dictionary)
 Let ${\cal H}$ be a complex Hilbert space and ${\cal A}$ a dictionary of ${\cal H}.$ Denote by
\[ \tilde{\cal A}=\cup_{k=0}^\infty {\cal A}_k\]
the \emph{Complete Dictionary Induced from} ${\cal A}.$  \\

Examples of such formulated complete dictionaries include the Complete Szeg\"o Dictionary $\tilde{\cal D}$ given in \S 1 (also see \cite{QWa}), where directional derivatives are replaced by the complex derivative, and one in (\cite{QSW}) for the quaternions case. We note that in the 1-D AFD case the orthonormal system functions $B_n$'s, no matter with repeated parameters or not, have an explicit and uniform formula in the parameters: the multiple indicates the order of the derivative. Due to the above analysis, from now on, to perform P-OGA we will be based on the complete dictionary $\tilde{\cal A}$ consisting of the elements of ${\cal A},$  as well as all possible first and higher order directional directives.

Now we show that P-OGA is applicable for the poly-disc Hardy $H^2$ space with the Complete Product-Szeg\"o Dictionary. We have the following

\begin{theorem}\label{P-OGA applicability}  For the poly-disc Hardy $H^2$ space with the Complete Product-Szeg\"o Dictionary $\tilde{{\cal D}^2}$ we can take $\rho =1$ in the Pre-Orthogonal $\rho$-Maximal Selection Principle (\ref{New}).\end{theorem}

\noindent {\bf Proof} Denote by $E_{\{a,b\}}$ a general element of $\widetilde{{\cal D}^2}.$ They are parameterized by
complex pairs $(a,b)$ in ${\bf D}^2.$  The elements in $\widetilde{{\cal D}^2}\setminus {\cal D}^2$ are not of the tensor product form, they, however, are finite linear combinations of tensor products of the form
\begin{eqnarray}\label{form} N(a,k)\partial_{\overline{a}}^ke_a \otimes N(b,l)\partial_{\overline{b}}^l e_b, \quad k+l>0,\end{eqnarray}
where $N(a,k)$ and $N(b,l)$ are the normalizing constants such that
\[  \| N(a,k)\partial_{\overline{a}}^k e_a \|=1, \quad  \| N(b,l)\partial_{\overline{b}}^l e_b \|=1.\]
Note that $N(a,0)=\sqrt{1-|a|^2}, N(b,0)=\sqrt{1-|b|^2}.$
Due to the singularity of $\partial_{\overline{a}}^k e_a$ when $|a|\to 1$ and that of $\partial_{\overline{b}}^l e_b$ when $|b|\to 1,$ we have
\begin{eqnarray}\label{limitzero}
\lim_{|a|\to 1} N(a,k)=0, \quad  \lim_{|b|\to 1} N(b,l)=0, \qquad k,l\ge 0.\end{eqnarray}
  It amounts to showing that there exists an interior pair $(a_n, b_n)\in {\bf D}^2$
such that
\begin{eqnarray}\label{availability}
& & \qquad  (a_n, b_n)=\arg \max \{ |\langle g_n, \tilde{B}_n \rangle | \ :\ (a,b)\in {\bf D}^2,
 \{B_1,...,B_{n-1}, \tilde{B}_n\} \nonumber \\ & & \ {\rm is\ the\ G-S\ orthogonalization \ of }\  E_{\{a_1,b_1\}},...,E_{\{a_{n-1},b_{n-1}\}}, E_{\{a,b\}}\},\end{eqnarray}
where $a_1,b_1, ..., a_{n-1},b_{n-1}$ are previously fixed. The $E_{\{a_n,b_n\}}$ will be the next selected element in the Complete Product-Szeg\"o Dictionary $\widetilde{{\cal D}^2}$ so to formulate $B_n.$ To prove the availability of (\ref{availability}) we first show
\begin{eqnarray}\label{next to show}  \lim_{|a|\to 1\ {\rm or}\  |b|\to 1} |\langle g_n, \tilde{B}_n \rangle |=0.\end{eqnarray}
Owing to the general relation (\ref{nonumber}) proved for $a\in \widetilde{{\cal D}^2}$ and the remark on the formulation of $E_{\{a,b\}},$ the above limit is a consequence of
\begin{eqnarray}\label{1st}  \lim_{|a|\to 1\ {\rm or}\  |b|\to 1} |\langle g_n, N(a,k)\partial_{\overline{a}}^ke_a\otimes N(b,l)\partial_{\overline{b}}^l e_b \rangle |=0\end{eqnarray}
and
\begin{eqnarray}\label{2nd}  \lim_{|a|\to 1\ {\rm or}\  |b|\to 1} \|Q_{n-1}(\partial_{N(a,k)\overline{a}}^ke_a \otimes N(b,l)\partial_{\overline{b}}^l e_b)\|=1,\end{eqnarray}
where $k+l\ge 0.$\\

Due to the orthogonality between
\[ N(a,k)\partial_{\overline{a}}^ke_a \otimes N(b,l)\partial_{\overline{b}}^l e_b \quad {\rm and }\quad
 {g}-\langle {g}, N(a,k)\partial_a^ke_a \otimes N(b,l)\partial_b^l e_b\rangle N(a,k)\partial_{\overline{a}}^ke_a \otimes N(b,l)\partial_{\overline{b}}^l e_b,\]
  we have
\[ \|g_n\|^2= |\langle g_n, N(a,k)\partial_{\overline{a}}^ke_a \otimes N(b,l)\partial_{\overline{b}}^l e_b \rangle|^2 + \|{g}-\langle {g}, N(a,k)\partial_{\overline{a}}^ke_a \otimes N(b,l)\partial_{\overline{b}}^l e_b\rangle N(a,k)\partial_{\overline{a}}^ke_a \otimes N(b,l)\partial_{\overline{b}}^l e_b\|^2.\]
Then (\ref{1st}) is a consequence of
\begin{eqnarray}\label{consequence}
 \lim_{|a|\to 1\ {\rm or}\  |b|\to 1}\|{g}-\langle {g}, N(a,k)\partial_{\overline{a}}^ke_a \otimes N(b,l)\partial_{\overline{b}}^l e_b\rangle N(a,k)\partial_{\overline{a}}^ke_a \otimes N(b,l)\partial_{\overline{b}}^l e_b\|^2=0.
\end{eqnarray}
Performing what is done in (\ref{inequality}), we are reduced to showing
 \begin{eqnarray}\label{reduced}   \lim_{|a|\to 1\ {\rm or}\  |b|\to 1}\|(P_r\otimes P_s) \ast (N(a,k)\partial_{\overline{a}}^ke_{a}\otimes N(b,l)\partial_{\overline{b}}^le_{b})\|\to 0.\end{eqnarray}
Now, with the fixed $0<r<1$ and $0<s<1$, since $N(a,k)\partial_{\overline{a}}^ke_a \otimes N(b,l)\partial_{\overline{b}}^l e_b\in H^2,$ there follows, for $z=re^{it}, w=se^{iu},$
\[ (P_r\otimes P_s) \ast (N(a,k)\partial_{\overline{a}}^ke_a \otimes N(b,l)\partial_{\overline{b}}^l e_b)(e^{it},e^{iu})=(N(a,k)\partial_{\overline{a}}^ke_{a})(z)(N(b,l)\partial_{\overline{b}}^le_{b})(w).\]
With explicit computation we have
\begin{eqnarray*}
\| (P_r\otimes P_s) \ast (N(a,k)\partial_{\overline{a}}^ke_a \otimes N(b,l)\partial_{\overline{b}}^l e_b)\|^2 = \frac{N(a,k)}{N(ra,k)}\frac{N(b,l)}{N(sb,l)}\to 0, \quad {\rm as} \quad |a|\to 1 \quad {\rm or}\quad |b|\to 1.\end{eqnarray*}
 With the relation (\ref{with}), the same reasoning gives
\[1=\| N(a,k)\partial_{\overline{a}}^ke_a \otimes N(b,l)\partial_{\overline{b}}^l e_b \|^2= \sum_{k=1}^{n-1} |\langle N(a,k)\partial_{\overline{a}}^ke_a \otimes N(b,l)\partial_{\overline{b}}^l e_b, B_k \rangle|^2 + \|Q_{n-1} (N(a,k)\partial_{\overline{a}}^ke_a \otimes N(b,l)\partial_{\overline{b}}^l e_b)\|^2.\]
By invoking (\ref{1st}),
\[ \lim_{|a|\to 1\ {\rm or}\  |b|\to 1}\sum_{k=1}^{n-1} |\langle N(a,k)\partial_{\overline{a}}^ke_a \otimes N(b,l)\partial_{\overline{b}}^l e_b, B_k \rangle|^2 =0.\]
We hence obtain (\ref{2nd}). Therefore, we have proven (\ref{next to show}). Next we indicate that we can always select an element of $\widetilde{{\cal D}^2},$ being a linear combination of finitely many
 directional derivatives of certain orders of the elements in ${\cal D}^2,$ that gives rise to the maximum of  $|\langle g_n, B^{a}_n\rangle |, a\in \widetilde{\cal A}$ (See the explanation before the definition of Induced Complete Dictionary).  The proof is complete.\\

The following theorem shows that in the one complex unit disc case the Complete Szeg\"o Dictionary defined in \S 1 through (\ref{def.ps}) coincides with the Complete Dictionary Induced by the Szeg\"o Dictionary defined in \S 4, and under which the related P-OGA reduces to 1-D AFD.\\

\begin{theorem}
In the 1-D unit disc context the P-OGA under the Complete Szeg\"o Dictionary $\tilde{\cal D}$ reduces to AFD. \end{theorem}

 \noindent{\bf Proof} In view of the relation
 \[\langle f_k, e_{a_k}\rangle =\langle g_k,B_k\rangle \]  cited in (\ref{reduce}) the AFD selection of $e_{a_k}$ in ${\cal D}$ in relation to the reduced remainder $f_k$ corresponds to the P-OGA selection of $B_k$ in relation to the orthogonal standard remainder $g_k.$  As indicated in \S 1, Remark 4, the system functions $B_k$ are formed through the G-S orthogonalization on elements in Complete Szeg\"o Dictionary,  we see that the P-OGA, in fact, is performed with respect to the Complete Szeg\"o Dictionary $\tilde{\cal D}.$  We note that a system function $B_k=B_{\{a_1,...,a_n\}}$ is parameterized directly by ${\bf a}=\{a_1,...,a_n\}\in \cup_{k=1}^\infty {\bf D}^k,$ where the multiple $m_k$ of $a_k, 1\leq k\leq n,$ indicates the order of the derivative of the Szeg\"o kernel that in use (See the definition of $E_k(z)$ in \S 1). The proof is complete.
\\

We finally note that all the results obtained in \S 2 and \S 3 for product spaces from several copies of the unit disc ${\bf D}$ have their counterparts in product spaces from several copies of the upper-half complex plane.
  \\

 \noindent{\bf Acknowledgement} The author would like to thank K-I. Kou, J-S. Huang, X-M. Li, L-M. Zhang, Q-H. Chen, W-X. Mai, L-H. Tan, W. Wu and Y. Gao for helpful discussions.

\end{document}